\documentclass[a4paper,reqno]{amsart}
\usepackage{amsmath}
\usepackage{amsfonts}
\usepackage{amsthm}
\usepackage{amssymb}
\usepackage{mathrsfs}
\usepackage{bm}
\usepackage{booktabs}
\usepackage{lipsum}
\usepackage{bussproofs}
\usepackage{hyperref}
\usepackage{tikz}
\usepackage{enumerate}
\usetikzlibrary{arrows,automata,positioning,decorations.markings,arrows.meta}
\usepackage{python}
\usepackage{comment}
\usepackage{amsmath}
\usepackage{amsfonts}
\usepackage{amsthm}
\usepackage{amssymb}
\usepackage{mathrsfs}

\def\parentheses#1{{\left( {#1} \right)}}
\def\of{\parentheses}

\def\Fraction#1/#2{\frac{#1}{#2}}
\def\fraction#1/#2{\tfrac{#1}{#2}}
\def\solidus#1/#2{%
    {{#1} \!\left/\!\vphantom{{#1}{#2}}\!\right. {#2}}%
}

\def\blackboard{\mathbb}

\def\integers{\blackboard Z}

\def\reals{\blackboard R}

\def\set#1{{\left\lbrace {#1} \right\rbrace}}
\def\setOf#1:#2{%
    \set{{#1} \,\left\vert\vphantom{{#1}{#2}}\right. {#2}}%
}

\def\emptySet{\varnothing}

\def\sequence#1{{\left\langle {#1} \right\rangle}}
\def\sequenceOf#1:#2{%
    \sequence{{#1} \,\left\vert\vphantom{{#1}{#2}}\right. {#2}}%
}

\def\familyOf#1:#2{{\parentheses{#1}_{#2}}}

\def\Union{\bigcup}
\def\intersect{\cap}

\def\setMinus{\setminus}

\def\functionSpace#1->#2{{{#2}^{#1}}}

\def\function#1:#2->#3{{{#1} \colon {#2} \to {#3}}}
\def\mapsTo{\mapsto}
\def\map#1->#2{\parentheses{{#1} \mapsTo {#2}}}
\def\into{\hookrightarrow}
\def\onto{\twoheadrightarrow}

\def\Restriction#1|#2{{\left.\! {#1} \right\vert_{#2}}} 

\def\presentation#1{{\left\langle {#1} \right\rangle}}
\def\presentationOf#1:#2{%
    \presentation{{#1} \,\left\vert\vphantom{{#1}{#2}}\right. {#2}}%
}

\def\normalClosure#1{{\left\langle\!\left\langle {#1} \right\rangle\!\right\rangle}}

\def\interval#1#2,#3#4{{\left#1 {#2}, {#3} \right#4}}

\def\modulus#1{{\left\lvert {#1} \right\rvert}}

\def\p{\parentheses}

\newtheorem{presentment}{}[section]

\theoremstyle{plain}
\newtheorem{theorem}[presentment]{Theorem}
\newtheorem{proposition}[presentment]{Proposition}
\newtheorem{lemma}[presentment]{Lemma}
\newtheorem{corollary}[presentment]{Corollary}

\theoremstyle{definition}
\newtheorem{definition}[presentment]{Definition}

\theoremstyle{remark}
\newtheorem*{claim}{Claim}
\newtheorem*{remark}{Remark}

\newtheorem*{acknowledgements}{Acknowledgements}

    \def\fraction#1/#2{\tfrac{#1}{#2}}
    \def\Fraction#1/#2{\frac{#1}{#2}}
\newcommand*{\BigFraction}{}
    \def\BigFraction#1/#2{\dfrac{\displaystyle{#1}}{\displaystyle{#2}}}

    \def\setOf#1:#2{\set{{#1} \,\left\vert\vphantom{{#1}{#2}}\right. {#2}}}

    \def\sequenceOf#1:#2{\sequence{{#1} \,\left\vert\vphantom{{#1}{#2}}\right. {#2}}}

    \def\function#1:#2->#3{{{#1} \colon {#2} \to {#3}}}

\tikzset{
  injective/.style={
    postaction={
      decorate,
      decoration={
        markings,
        mark=at position 1 with {\arrow{>}},
        mark=at position 0 with {\node[inner sep=0pt] {$\hookrightarrow$};}
      }
    }
  }
}

\theoremstyle{plain}
\newtheorem{lettertheorem}{Theorem}

\theoremstyle{plain}
\newtheorem{lettercorollary}[lettertheorem]{Corollary}

\usepackage{pst-node}
\usepackage{tikz-cd} 
\tikzset{node distance=2cm, auto}

\allowdisplaybreaks[0]
\title{One-Relator Products of Left-Orderable Groups} 

\author
{Alba González}
\address
{Mathematical Institute, University of Oxford, Andrew Wiles Building,
Radcliffe Observatory Quarter,
Woodstock Road,
Oxford,
OX2 6GG}
\email {} \email {alba.gonzalez@maths.ox.ac.uk}

\author
{Tejas Mittal}
\address
{Mathematical Institute, University of Oxford, Andrew Wiles Building,
Radcliffe Observatory Quarter,
Woodstock Road,
Oxford,
OX2 6GG}
\email {} \email {tejas.mittal@maths.ox.ac.uk}

\author
{Andrew Wood}
\address
{Mathematical Institute, University of Oxford, Andrew Wiles Building,
Radcliffe Observatory Quarter,
Woodstock Road,
Oxford,
OX2 6GG}
\email {} \email {andrew.wood@maths.ox.ac.uk}

\begin{document}

\begin{abstract}
We prove that, for left-orderable groups $A,B$ and a word $w\in A*B$ which is not conjugate into $A$ or $B$ and is not a proper power, the one-relator product $A*B/\normalClosure{w}$ is torsion-free, answering a question posed by Chen and Lodha. We do this by proving a relative asphericity result for the one-relator product. We also generalise the B. B. Newman Spelling Theorem to one-relator products and deduce that a one-relator product $A*B/\normalClosure{w^k}$ of left-orderable groups with torsion is hyperbolic relative to $\{A,B\}$.
\end{abstract}

\maketitle
\section{Introduction}

Recently, Chen and Lodha resolved the Wiegold problem from $1976$ in the affirmative, which asked whether there exists finitely generated perfect groups with normal rank strictly greater than one~\cite{wiegold}. Their examples come from free products of left-orderable groups. Recall that a group $G$ is \emph{left-orderable} if there is a total strict order $<$ on $G$ such that $f g < f h$ for each $f \in G$ whenever $g < h$. More precisely, they proved the following theorem.

\begin{theorem}[{\cite[Theorem~B]{wiegold}}]
\label{theorem:wiegold}
For any non-trivial left-orderable groups $A$ and $B$ and any $w\in A\ast B$ not conjugate into $A$, the natural map ${A\to A\ast B/\normalClosure{w}}$ is injective. In particular, free products of non-trivial left-orderable groups have normal rank at least two.
\end{theorem}

A natural follow-up question posed by Chen and Lodha is whether the one-relator product $A \ast B / \normalClosure{w}$ is torsion-free, where $w \in A \ast B$ is not conjugate into $A$ or $B$ and is not a proper power~\cite[Question $6.2$]{wiegold}. We prove that the relative presentation $A*B/\normalClosure{w}$ is aspherical and deduce that the one-relator product $A*B/\normalClosure{w}$ must indeed be torsion-free.

\begin{lettertheorem}\label{lettertheorem:asphericity}
Let $A, B$ be left-orderable groups and $w\in A\ast B$ be an element which is not conjugate into $A$ or $B$ and is not a proper power. Then, the relative presentation $A\ast B/\normalClosure{w}$ is aspherical.
\end{lettertheorem}

The heart of our proof involves an analysis of oriented pictures. Such pictures are also the main tool in an alternative proof of Theorem \ref{theorem:wiegold} discovered independently by Hamish Short (private communication).

\begin{lettercorollary}[Corollary \ref{corollary:torsionfree}]
\label{lettercorollary:torsionfree}
 Let $A, B$ be left-orderable groups and $w\in A\ast B$ be an element which is not conjugate into $A$ or $B$ and is not a proper power. Then, the one-relator product $A\ast B/\normalClosure{w}$ is torsion-free.
\end{lettercorollary}

While proving Theorem \ref{lettertheorem:asphericity}, we will construct an explicit classifying space for the one-relator product $A\ast B/\normalClosure{w}$ (see Theorem \ref{theorem:asphericity}). Following an argument of Howie~\cite{Howie1984CohomologyOO}, this immediately allows us to compute (co)-homology for the one-relator product. The proof for one-relator products of left-orderable groups is identical to the analogous result for one-relator products of locally indicable groups. We refer the reader to~\cite{Howie1984CohomologyOO} for the proof and a more precise version of the statement.

\begin{lettercorollary}[Corollary \ref{corollary:cohomologycompute}]
\label{lettercorollary:cohomologycompute}
 Let $A, B$ be left-orderable groups and $w\in A\ast B$ be an element which is not conjugate into $A$ or $B$ and is not a proper power. Let $G$ be the one-relator product $A\ast B/\normalClosure{w}$ and $\mathbf{1}$ denote the trivial subgroup of $G$. Then, for $k\geq3$,
 \[H^{k}(G;\text{--}) = H^{k}(A;\text{--}) \oplus H^{k}(B;\text{--}), \]
 \[H_k(G;\text{--}) = H_k(A;\text{--}) \oplus H_k(B;\text{--}).\]
 For $k = 1,2$, the (co)-homology groups of $G$ can be computed via the following exact sequences:
\[0\rightarrow H^{1}(G;\text{--}) \rightarrow H^{1}(A\ast B; \text{--})\rightarrow H^{0}(\mathbf{1}; \text{--})\rightarrow H^{2}(G; \text{--})\rightarrow H^{2}(A; \text{--})\oplus H^{2}(B; \text{--})\rightarrow 0,\]
\[0\rightarrow H_{2}(A; \text{--})\oplus H_{2}(B; \text{--}) \rightarrow H_{2}(G; \text{--})\rightarrow H_{0}(\mathbf{1}; \text{--})\rightarrow H_{1}(A\ast B; \text{--})\rightarrow H_{1}(G; \text{--})\rightarrow 0.\]
\end{lettercorollary}

We also generalise some classical results in one-relator group theory to one-relator products of left-orderable groups. Firstly, we prove a version of Weinbaum's Subword Theorem stating that proper, non-empty, alternating subwords of $w$ do not represent the identity in $A\ast B/\normalClosure{w}$. This generalises Corollary 6.1 in \cite{wiegold}.

\begin{lettertheorem}[Theorem \ref{theorem:subword}]
\label{lettertheorem:subword}
Let $A, B$ be left-orderable groups and $w\in A\ast B$ be a word which is not conjugate into $A$ or $B$. Then, for any proper, non-trivial subword $u$ of $w$, we have $u \notin \normalClosure{w}$.
\end{lettertheorem}

Secondly, we establish a version of the B. B. Newman Spelling Theorem for one-relator products of left-orderable groups.  The analogous result for one-relator products of locally indicable groups was proved by Howie and Pride \cite{HowiePride}. We remark that our proof strategy using oriented pictures is different from their's, and it is distinct from classical approaches for one-relator groups. As a corollary, we get that one-relator products of left-orderable groups with torsion are relatively hyperbolic.

\begin{lettertheorem}[Theorem \ref{theorem:newmanspellingthm}]
\label{lettertheorem:spelling}
Let $A$, $B$ be left-orderable groups and $k>1$ an integer. Consider the one-relator product $G=A*B/\normalClosure{w^k}$ where $w$ is not a proper power and not conjugate into $A$ or $B$. Let $u$ be a reduced non-trivial word in $A \ast B$ representing the identity element in $G$. Then, $u$ contains an alternating subword $v$ so that either $v$ or $v^{-1}$ is a subword of $w^k$ with 
\[
\emph{length}\of{v} > \Fraction k - 1 / {k} \emph{length}\of{w^k}.
\]
\end{lettertheorem}

\renewcommand{\thepresentment}{\Alph{presentment}}

\begin{lettercorollary}[Corollary \ref{relativehyperbolicity}]
\label{corollary:relativelyhyperbolic}
 For left-orderable groups $A,B$, a one-relator product $A*B/\normalClosure{w^k}$ with torsion is hyperbolic relative to $\{A,B\}$.
\end{lettercorollary}

\begin{remark}
 Right-orderability is defined by requiring that the total order is invariant under right multiplication. Note that a group is left-orderable if and only if it is right-orderable. For the rest of the paper, we will use right-orderability as it aligns with Chen and Lodha's work in \cite{wiegold}.
\end{remark}

\noindent \textit{Organisation of the paper.} In Section~\ref{section:orientedpictures}, we introduce oriented pictures, which we use to give an alternative proof of {\cite[Theorem~B]{wiegold}}. We also include some technical results needed to prove our theorems. Section~\ref{section:asphericity} is devoted to proving Theorem \ref{theorem:asphericity} and some of its consequences, including Corollary \ref{lettercorollary:torsionfree}. In section \ref{section:classical}, we prove Theorem \ref{lettertheorem:subword}, Theorem \ref{lettertheorem:spelling}, and Corollary \ref{corollary:relativelyhyperbolic}.
\\

\noindent \textit{AI statement}: No AI tools were used to develop the ideas, arguments, figures, results or proofs
in this article, nor to draft the manuscript. AI (ChatGPT 5.6) was used for locating references. The authors take full responsibility for the results.

\begin{acknowledgements}
The authors are very grateful to Martin Bridson for carefully proofreading this manuscript, and for the many helpful conversations and advice on this work. The authors are also in debt to Hamish Short, Henry Wilton, Lvzhou Chen, and Marco Linton for comments on a previous version of this paper. In particular, we thank Short for emphasising the importance of oriented pictures, leading to a clearer exposition of the paper; Wilton for fruitful discussions on stacking and for his mini-course in the conference about Aspherical 2-complexes at ICMAT, Madrid; Chen for fruitful discussions during the aforementioned conference; and Linton for suggesting that we think about computing (co)-homology groups for the one-relator product. Finally, the authors are thankful to Ojas Mittal for his help with a concise argument for Lemma \ref{homotopy} and Colin McCulloch for helpful conversations and pointing out a subtlety in the proof of Theorem $\ref{theorem:asphericity}$.
\end{acknowledgements}

\renewcommand{\thepresentment}{\thesection.\arabic{presentment}}

\section{Oriented Pictures and background}\label{section:orientedpictures}
In this section,  we give the necessary background for what is required to understand the paper. We give a proof of Chen and Lodha's main result~\cite[Theorem~B]{wiegold} in the language of pictures with oriented edges. These oriented pictures will be the key to proving the results in this paper.  We now give our definition of pictures (see~\cite{pictures, warwick177864, Chinyere_2016} for further background on the classical definition of pictures). This variation on the classical notion of pictures came into prominence due to Chen and Lodha \cite{wiegold}, while the idea of pictures with oriented edges goes back to Duncan and Howie \cite{DuncanHowie1991}.

First, throughout we let $G = A \ast B$ where $A$ and $B$ are non-trivial groups, and let $w$ be a cyclically reduced non-trivial element of $G$ which is not conjugate into $A$ or $B$. Let $K$ be the cell complex obtained from classifying spaces $K\of{A, 1}$ and $K\of{B, 1}$ with base points $p_A$ and $p_B$ by adding a 1-cell $e$ connecting $p_A$ to $p_B$. We take the midpoint of $e$ to be our basepoint. Call this basepoint $\ast$.

\begin{definition}\label{definition:picture}
We call $\Gamma = \p{f, S}$, where $S$ is a compact oriented surface and $f$ is a map from $S$ to $K$,  a \emph{picture on $S$} (relative to $\set{A, B}$) if the image of all except possibly one boundary component of $S$ represents the conjugacy class of $w^{\pm 1}$. A component of this type is called a \emph{$w$-circle}, and the union of $w$-circles is called the \emph{$w$-boundary}. The other type of boundary (if it exists) represents the conjugacy class of an arbitrary element in $G$, and we call this the \emph{boundary} of $\Gamma$, denoted by $\partial \Gamma$.  The \emph{area} of $\Gamma$ is the number of $w$-circles in $\Gamma$, denoted by $\text{area}\of{\Gamma}$. In the case that $w$ is a proper power, say $w = u^k$ (where $u$ is not a proper power), we define the \emph{degree} of $\Gamma$ to be $\modulus{k} \text{area}\of{\Gamma}$ and denote this by $\deg\of{\Gamma}$. 
\end{definition}

We remark that Definition~\ref{definition:picture} is similar to the definition of $w$-admissible surfaces (see~\cite[Definition $3.1$]{wiegold}). More precisely, a picture with boundary in $A$ (or $B$), and with non-empty $w$-boundary is a $w$-admissible surface without boundary components representing proper powers of $w^{\pm 1}$, and with a single $A$-boundary (or $B$-boundary) component.

For the definition of area/degree to be useful, we need to consider the following definition.

\begin{definition}\label{definition:compressibility}
We call a picture $\Gamma = (f, S)$ \emph{compressible} if there is an embedded subsurface $P\subset S$ homeomorphic to a pair of pants, such that two boundary components of $P$ are on the $w$-boundary of $S$ and represent the conjugacy class of $w$ and $w^{-1}$ respectively, and the third boundary component of $P$ is in the interior of $S$ and represents a loop homotopic to the identity under the map $f$. We say $\Gamma$ is \emph{incompressible} if there is no such $P$.
\end{definition}

We remark that for any group element $u \in \normalClosure{w}$, we get a planar incompressible picture $\Gamma$ with $\partial\Gamma$ representing $u$. This is the dual of the relative Van Kampen diagram for $u$~\cite[Part III: Diagrams and groups]{pictures}. 

Let $\Gamma = \p{f, S}$ be a picture. Similar to~\cite[Section~2]{wiegold}, we can assume that the pre-image $f^{-1}(\ast)$ is a disjoint union of arcs using transversality.  This equips our pictures $\Gamma$ with arcs between boundary components.\\ 

Let $L$ be any CW-complex obtained by attaching a number of cells to $K$ so that the loop $l$ corresponding to the word $w$ bounds some disc $D_l$. Then, for any picture $\Gamma = (f,S)$, we can extend $f$ to a map $\overline{f}:\overline{S}\to L$ where $\overline{S}$ is now a surface (with one boundary $\partial \Gamma$) where we glue discs to each $w$-circle. These \textit{$w$-discs} are then mapped to the disc $D_l$ via $\overline{f}$. We note that compressing a picture $\Gamma = (f,S)$ to another one $\Gamma' = (f', S')$, by deleting the corresponding pair of pants, induces a homotopy equivalence between the maps $\overline{f}:\overline{S}\to L$ and $\overline{f'}:\overline{S'}\to L$. This is analogous to reductions in the classical notion of pictures and the result below is folklore, but as we could not find a reference in the literature, we include its proof, for the sake of completeness.

\begin{lemma}
\label{homotopy}
Let $\Gamma = (f, S)$ be a compressible picture with an embedded pair of pants $P$ and let $\Gamma'=(f',S')$ be the picture obtained by deleting $P$ from $S$. Then, the maps $\overline{f'}$ and $\overline{f}$ are homotopic.
\end{lemma}

\begin{proof}

Let $P$ bound the loops $l_1$, $l_2$, and $l_3$ with $l_1$, $l_2$ representing free homotopy classes corresponding to $w$ or $w^{-1}$ and $l_3$ representing the free homotopy class of the identity under $f$.

Thus, up to a homotopy, we may assume $l_3$ maps to the point $\ast$ under $\overline{f}$. Let $l_3$ bound a disc $D$ in $\overline{S}$.

\begin{figure}[ht]
    \centering
    \includegraphics[width=100mm]{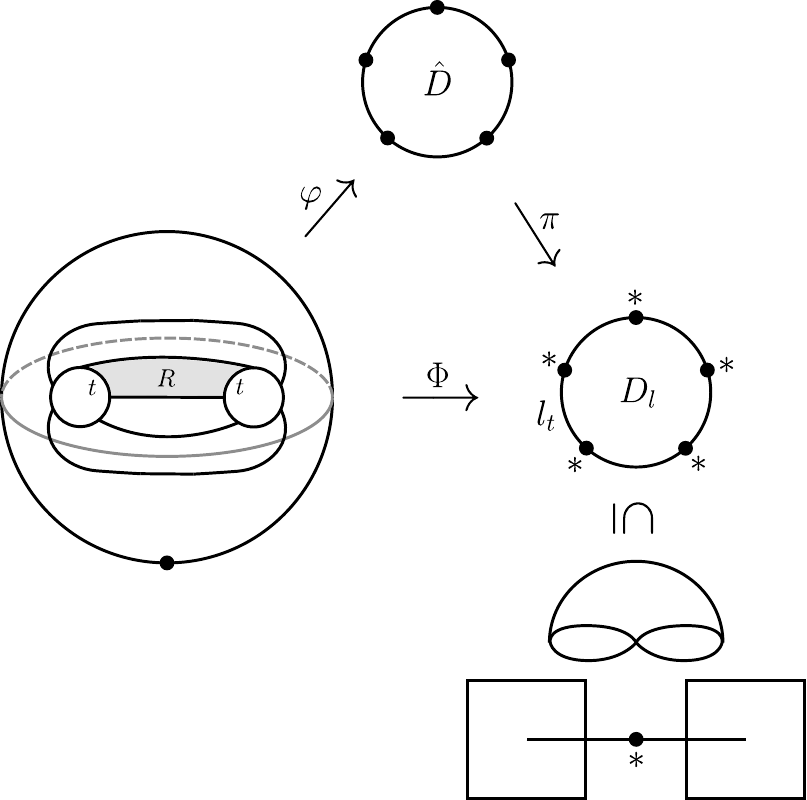}  
    \caption{\label{fig:Homotopy}The map $\Phi$ is homotopic to the constant map $c_\ast$ since it factors through a contractible space (a disc).}
\end{figure}

We now show that $\overline{f}\big|_D \simeq \operatorname{c}_{\ast} \text{ rel } \partial D$, where $\operatorname{c}_{\ast}$ sends every point to $\ast$. Make $f\big |_D$ transverse to $\ast$ so that $f\big|_D^{-1} (\ast)$ is a 1-manifold. The pre-image $f\big|_D^{-1} (\ast)$ might have some embedded loops in the interior $\mathring{D}$ of $D$, but these can be removed up to a homotopy $\text{ rel } \partial D$ as follows. 

Consider an inner-most loop $l$ which bound a disc $D'\subset D$. If $D'$ contains neither of the loops $l_1$ or $l_2$, then $\overline{f}(D')\subseteq K$. Since $K$ is aspherical, $\overline{f}\big|_{D'}$ is homotopic $\text{ rel } \partial D'$ to the constant map sending $D'$ to the point $\ast$. Via a homotopy, shrink $D'$ to a point $q$ which gets mapped to $\ast$. Moreover, for small $\varepsilon$, $N_{\varepsilon}(q)\setminus q$ is connected and its image does not contain $\ast$. So, $\overline{f}(N_{\varepsilon}(q))$ lies in one of $K(A,1)$ or $K(B,1)$ union a half-edge and can be pushed away from $\ast$ via a homotopy. This removes the loop $l$.\\
If $D'$ contains just one of the loops among $l_1$ or $l_2$, then the loop $\overline{f}(l_i)$ is freely homotopic to $\overline{f}(\partial D')=\ast$. However, $\overline{f}(l_i)$ represents $w$ or $w^{-1}$ which is not the identity in $A*B$; this is a contradiction.\\
If $D'$ contains both of the loops $l_1$ or $l_2$, then we are in the case when $f\big |_{\mathring{D'}}^{-1} (\ast)$ has no embedded loops which we will handle now.

Now, $\overline{f}\big |_{\mathring{D}}^{-1} (\ast)$ is a union of disjoint arcs. Moreover, these arcs must start and end at $l_1$ or $l_2$. Thus, the arcs must be as in Figure \ref{fig:Homotopy}, otherwise some arc would start and end at the same loop and we obtain a sub-word of $w$ representing the trivial element in $A*B$. Let $\Phi$ be the map corresponding to $\overline{f}\big |_D$ obtained by collapsing $\partial D$ to a point. Thus, $\Phi$ is a map from a sphere $S'$ to $L$ and proving $\overline{f}\big|_D \simeq \operatorname{c}_{\ast} \text{ rel } \partial D$ amounts to showing $\Phi$ is homotopic to the constant map. Let $R$ be a region bounded by two of the arcs and segments representing the same letter $t$ of $w$ (see Figure \ref{fig:Homotopy}). Since $K$ is aspherical, one can re-define the map $\Phi$ on each region $R$ relative to $\partial R$ without changing the homotopy type of $\Phi$. Note that $\Phi(\partial R)=l_t$, where $l_t$ is the path in $L$ tracing the letter $t$ of $w$. Therefore, upto homotopy, we can assume that the inside of the region $R$ is also mapped inside $l_t$.

Finally, let $\varphi$ send every region $R$ to the corresponding path $l_t$ in $\partial\hat{D}$ and the two $w$-cells to the disc $\hat{D}$. Let $\pi\colon \hat{D} \to L$ corresponds to the attaching map for disc $D_l$.
Note that $\Phi$ factors through the disc $\hat{D}$ , i.e. $\Phi = \pi \circ \varphi$ as in Figure \ref{fig:Homotopy}. Since the disc $\hat{D}$ is contractible, $\varphi$ is homotopic to a constant map. Therefore, $\Phi = \pi \circ \varphi$ is homotopic to a constant map as well.

Finally, note that $\overline{f'}$ is the map obtained by collapsing the disc $D$ to a point. By the claim above, $\overline{f}\big|_D \simeq \operatorname{c}_{\ast} \text{ rel } \partial D$ and $\overline{f} = \overline{f'}$ on $S\setminus D$. By glueing the two homotopies, we get that $\overline{f} \simeq \overline{f'}$.
\end{proof}

Assigning orientations to the arcs between $w$-circles yields the aforementioned oriented pictures. This is done using relative stackings (as in~\cite[Section $5$]{wiegold}), which in turn was inspired by work of Louder and Wilton on stackings \cite{Louder_Wilton_2017}. In fact, as mentioned by Wilton at the conference on Aspherical $2$-complexes at ICMAT (Madrid, $2026$), using stackings one can give alternate geometric proofs to classical results in one-relator group theory (including, for example, the B.B. Newman's Spelling Theorem).   Chen and Lodha prove the following.

\begin{theorem}[{\cite[Theorem~F]{wiegold}}]\label{theorem:relative-stacking-exist}
Let $G = A \ast B$ be a free product of countable right-orderable groups $A$ and $B$.  For any $w 
\in G$ of cyclically reduced length at least $2$, a relative stacking of $w$ exists if and only if $w$ is not a proper power.
\end{theorem}
  
We now define what is meant by a relative stacking. Firstly, consider an action $\rho \colon G = A \ast B \to \text{Homeo}^+\of{\reals}$.  For a cyclically reduced word $w = a_1 b_1 \ldots a_n b_n$ with each $a_i \in A \setMinus \set{1_A}$ and $b_i \in B \setMinus \set{1_B}$, the \emph{trajectory} of a given $x \in \reals$ under $w$, is the set:
\[
\Omega\of{w, x} := \setOf{x \cdot \rho\of{u}}:{\text{$u$ is a non-empty prefix of $w$}},
\]
where a non-empty prefix $u$ of $w$ is of the form $u = a_1 b_1 \ldots a_i b_i$ or $u = a_1 b_1 \ldots a_j$, for $i, j \in \set{1, \ldots, n}$. We say the trajectory $\Omega\of{w, x}$ of $x$ under $w$ is \emph{stable}, if 
\begin{enumerate}
\item each element of $\Omega\of{w, x}$ occurs precisely once (i.e., the action of non-empty prefixes of $w$ on $x$ are unique); 
\item $x \cdot \rho\of{w} = x$. 
\end{enumerate}
A \emph{relative stacking} of $w$ is a right action $\rho \colon G = A \ast B \to \text{Homeo}^+\of{\reals}$ such that there exists $x \in \reals$ with a stable trajectory under $w$.

\begin{figure}[ht]
    \centering
    \includegraphics[width=90mm]{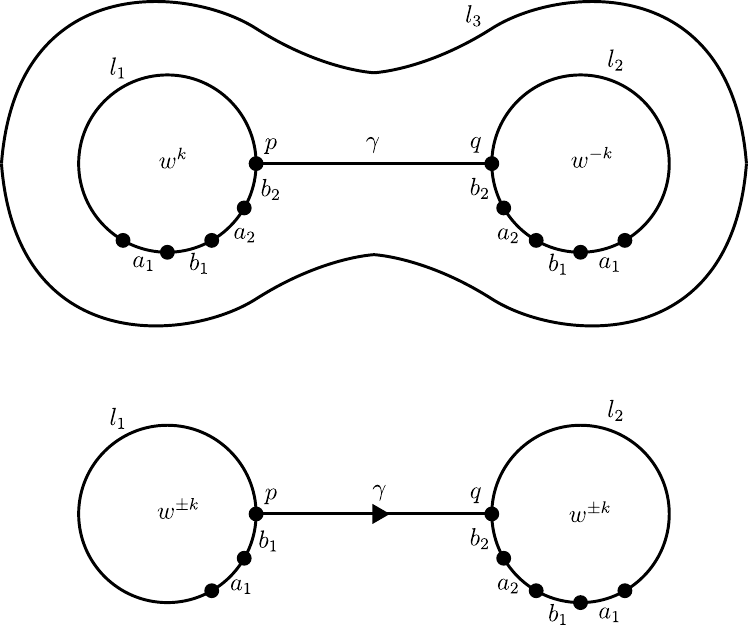} 
    \caption{\label{fig:OrientingArcs}The top part of the figure shows that if $u_p = u_q$, then $\Gamma$ is compressible. The bottom part tells how to orient the arc $\gamma$ using the relative stacking, when $x \cdot \rho\of{u_p} > x \cdot \rho\of{u_q}$.}
\end{figure}

Let $G = A \ast B$, where $A$ and $B$ are right-orderable groups, and let $w \in G$ such that $w$ is cyclically reduced of length at least $2$, and is not a proper power. We now explain how to orient a picture $\Gamma = (f,S)$ relative to $\set{A,B}$ with $w^k$-circles, where $k\geq 1$. In other words, we are looking at the one-relator product $A\ast B / \normalClosure{w^k}$.

By Theorem~\ref{theorem:relative-stacking-exist}, there exists a relative stacking $\rho$ of $w$;  let $x \in \reals$ be such that it has a stable trajectory under $w$. Via transversality, we may assume that the pre-image $f^{-1}(\ast)$ is a disjoint union of arcs (see~\cite[Section~2]{wiegold}). We are able to orient arcs joining $w^k$-circles as follows. Suppose we have an arc $\gamma$ connecting points $p$ and $q$ on $w^k$-circles $\ell_1$ and $\ell_2$, respectively (see Figure~\ref{fig:OrientingArcs}). Let $u_p$ and $u_q$ be the prefixes of $w$ corresponding to $p$ and $q$, respectively. More precisely, if $p$ is separating two segments reading $w$, then we define $u_p = e$. Otherwise, there will be a unique segment reading $w$ which contains $p$, and we define $u_p$ to be the prefix of $w$ read when following $\ell_1$ along this segment until $p$.

We orient the arc $\gamma$ from $p$ to $q$ if, and only if, $x \cdot \rho\of{u_p} > x \cdot \rho\of{u_q}$.  To see this orientation is well-defined, we must consider the case $x \cdot \rho\of{u_p} = x \cdot \rho\of{u_q}$.  Since $\rho$ is a relative stacking, by definition this forces $u_p = u_q$.  By applying the argument presented in~\cite[Lemma $4.3$]{wiegold}, this means $S$ is compressible (as witnessed by the loop $\ell_3$ in Figure~\ref{fig:OrientingArcs}), which is a contradiction.

We now give an alternative proof to \cite[Theorem~B]{wiegold} using oriented (planar) pictures and the lemma about oriented planar graphs below.

\begin{definition}
     Let $\Gamma$ be a finite, planar/spherical, oriented graph with possibly multiple edges between every pair of vertices. We say that a vertex is a \textit{sink} if all its edges are oriented inwards. Similarly, we say that a vertex is a \textit{source} if all its edges are oriented outwards. We call a cycle in $\Gamma$ an \textit{elementary cycle} if it corresponds to the boundary of a face of $\Gamma$.
\end{definition}

\begin{lemma}
\label{elementarycycle}
Let $\Gamma$ be a finite, oriented graph on a sphere or a plane with at least one edge.
\begin{enumerate}
    \item If at most one vertex of $\Gamma$ is a source or a sink, then $\Gamma$ contains an elementary cycle.
    \item If $\Gamma$ is on a sphere $S$ and no vertex of $\Gamma$ is a source or a sink, then $\Gamma$ contains at least two elementary cycles.
\end{enumerate}
\end{lemma}
\begin{proof}

     If there is a source/sink in $\Gamma$, take $v$ to be such vertex. Without loss of generality, assume $v$ is a source. Now, as all other vertices in $\Gamma$ admit both incoming and outgoing edges, by starting at a vertex $v'\neq v$ and following a path of outgoing edges, we must eventually repeat a vertex and hence get an oriented cycle. Suppose $C$ is a cycle which bounds a disc $D\subseteq S\setminus \set{v}$ containing minimal number of faces of $\Gamma$. If $D$ contains no vertex in its interior, then we have found our elementary cycle. Otherwise, since $v\notin D$, we start at a vertex $u$ inside $D$ and again follow a path of outgoing edges to get a cycle bounding a smaller number of faces. This contradicts the minimality of $C$.
     
     For (2), first get an oriented cycle in $S$ as above and then pass to a shortest such cycle $C$, so that it is embedded. Now, $C$ bounds two discs $D_1$ and $D_2$, one on each side of $C$ in the sphere $S$. Each of $D_1$ and $D_2$ has the structure of a finite planar orientated graph, all of whose vertices are not sources or sinks, and hence, by (1), must contain elementary cycles.
\end{proof}

We distil the arguments in \cite[Lemma 4.5 and Lemma 4.6]{wiegold} into the following proposition establishing that the oriented pictures we are concerned with do not admit sources/sinks or elementary cycles.

\begin{proposition}
\label{proposition:nosourcesinkcycle}
    Let $A$ and $B$ be non-trivial right-orderable groups, and let $w$ be an element in $G=A\ast B$ not conjugate into $A$ and not a proper power. Any incompressible oriented picture over $A\ast B/\normalClosure{w^k}$, thought of as a graph, contains no source/sink and no elementary cycle.
\end{proposition}

\begin{proof}
A source (resp. sink) would correspond to all arcs from a $w^k$-circle being oriented outwards (resp. inwards). However, the orientations are given via the relative stacking $\rho$ and $x \cdot \rho(u)$ must take minimal and maximal values for exactly two of the prefixes, which must then have different orientations.

Now suppose we have an elementary cycle $C$ bounding a disc piece $P$ whose polygonal boundary alternates between arcs and segments: $\gamma_1, s_1, \gamma_2, s_2,...,\gamma_m, s_m$ (see Figure \ref{fig:ElementaryCycle}).

\begin{figure}[ht]
    \centering
    \includegraphics[width=60mm]{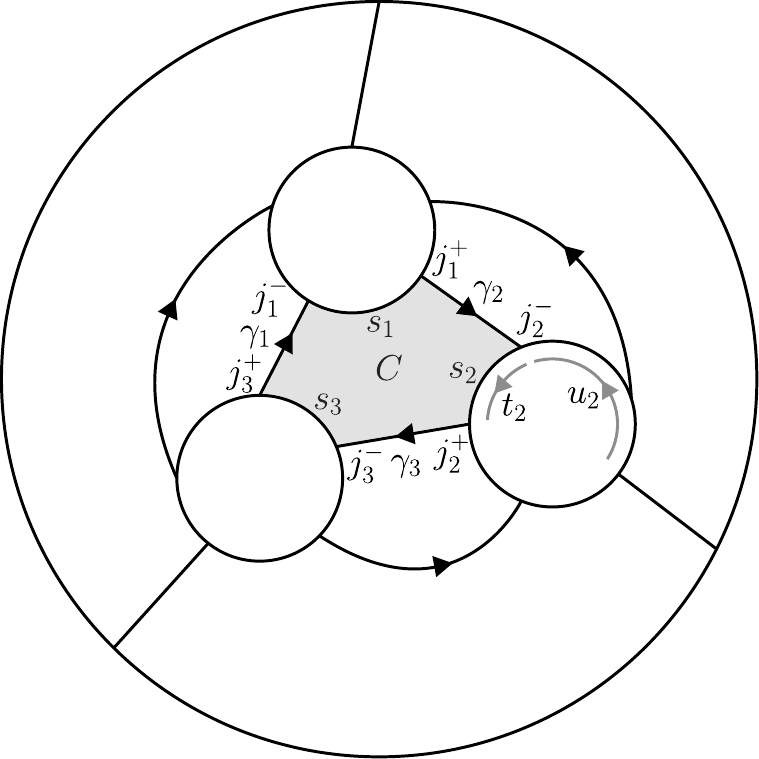}  
    \caption{\label{fig:ElementaryCycle}Pictures oriented via a relative stacking cannot have elementary cycles.}
\end{figure}

Indeed, suppose the segment $s_i$ starts and ends at junctures $j^{-}_i$ and $j^{+}_i$ which correspond to prefixes $u_i$ and $u_i t_i$ respectively. Without loss of generality, assume $\gamma_i$ is oriented from $j^{+}_i$ to $j^{-}_{i+1}$ so that we have $x\cdot \rho(u_i t_i)>x\cdot \rho(u_{i+1})$. Finally, note that $x\cdot \rho(u_1 t_1)>x\cdot \rho(u_2)$ and hence $x\cdot \rho(u_1 t_1 t_2) = x\cdot \rho(u_1 t_1)\rho(t_2)> x\cdot \rho(u_2) \rho(t_2) > x\cdot \rho(u_3)$. Applying this inductively, we have that for all $l\leq m$ (where we define $u_{m+1}$ to be $u_1$), $$x\cdot \rho(u_1 t_1 t_2 ... t_l) > x\cdot \rho(u_l) \rho(t_l)>x\cdot \rho(u_{l+1}).$$

Setting $l=m$ gives $x\cdot \rho(u_1 t_1 t_2 ... t_m) > x\cdot \rho(u_1)$. However, as $P$ bounds a disc, we have $t_1t_2...t_m = 1$ in $A \ast B$. Therefore, $x\cdot \rho(u_1) = x\cdot \rho(u_1 t_1 t_2 ... t_m)>x\cdot \rho(u_1)$. This contradiction implies that there are no elementary cycles.
\end{proof}

\begin{proof}[\textbf{Alternative proof of {\cite[Theorem~B]{wiegold}}}]
Using our notation, we will prove $A$ injects into $A\ast B/\normalClosure{w^k}$ where $w$ is not a proper power. Suppose, $a\in \normalClosure{w^k}$. This gives us an incompressible oriented picture $\Gamma$ with $\partial \Gamma$ reading $a$. Proposition \ref{proposition:nosourcesinkcycle} ensures that $\Gamma$ will not contain any source/sinks or elementary cycles. Since this picture witnesses $a\in \normalClosure{w^k}$ and $a\in A$, the only arcs will be amongst the $w^k$-circles and the picture can be thought of as an oriented planar graph. This contradicts Lemma \ref{elementarycycle}.
\end{proof}

Chen and Lodha prove Theorem \ref{theorem:wiegold} by establishing the following bound on the complexity of a $w$-admissible surface (see \cite[Definition $3.1$]{wiegold}). 

\begin{theorem}[{\cite[Theorem~D]{wiegold}}]
\label{thmA: characteristic estimate}
Let $G=A\ast B$ be a free product of right-orderable groups $A$ and $B$, and let $w$ be an element of $G$ not conjugate into $A$ or $B$. Then for any boundary-incompressible $w$-admissible surface $S$ without $2$-sphere or disk components, we have
    $$-\chi(S)\geq \deg(S).$$
\end{theorem}

 We now adapt the above Theorem~\ref{thmA: characteristic estimate} to oriented pictures. Recall that this means that one of the boundary components of the corresponding surface is allowed to be labelled $u$. As a consequence, one can prove relative hyperbolicity for one-relator products of right-orderable groups with torsion (see the \hyperref[Remark:relativehyperbolicity]{remark} below Corollary~\ref{relativehyperbolicity}). The following result is proved for locally indicable groups in~\cite[Theorem $3.3$]{DuncanHowie1991}.
 
Let $u \in A \ast B$ be such that 
\[
u = c_1 d_1 \cdots c_m d_m, \quad \text{where } c_i \in A \setMinus \set{1_A}, d_i \in B \setMinus \set{1_B} \text{ for every $i \in \set{1, \ldots, m}$}. 
\]
Then the \emph{length of $u$}, denoted by $\text{length}\of{u}$, is $2m$.

\begin{theorem}
\label{theorem:characteristicestimateboundary}
Let $G=A\ast B$ be a free product of right-orderable groups $A$ and $B$, and let $w$ be an element of $G$ not conjugate into $A$ or $B$ and not a proper power. Let $\Gamma = (f,S)$ be a picture over $A\ast B/\normalClosure{w^k}$ with $\partial \Gamma$ reading the word $u = c_1 d_1 \cdots c_m d_m \in A*B$ with $c_i\in A \setminus \set{1_A}$, $d_i\in B \setminus \set{1_B}$. Then,
\[
-\chi\of{S} \geq \deg\of{S} - \Fraction {1} / {2} \emph{length}\of{u}.
\]
\end{theorem}

\begin{proof}
Our proof strategy is similar to~\cite[Theorem $3.3$]{DuncanHowie1991} and~\cite[Theorem~D]{wiegold}.   Let $f \colon S \to K$ be the picture. As explained above, via transversality, we may assume that the pre-image $f^{-1}(\ast)$ is a disjoint union of arcs (see~\cite[Section~2]{wiegold}), and by using the relative stacking, we orient the arcs connecting any two $w^k$-circles.

Cutting the surface along the arcs gives us a disjoint union of pieces. Each piece $P$ is either a disc or an annulus with $\partial \Gamma$ on the inside reading a letter in $A$ or $B$. By Mayer-Vietoris, we have $\chi(S)=\sum_{\text{all pieces}} \chi(P)-E$, where $E$ is the number of arcs in $S$. Note that every arc contains exactly two junctures, so that we have 
\[
\chi(S)=\sum_{\text{all pieces}} \chi(P)-\Fraction {1} / {2} (\text{length}\of{u} + \modulus{w} \deg(S)).
\]
Call a piece an \emph{inner piece} if all the arcs on its boundary are oriented, i.e. the corresponding arcs connect two $w^k$-circles. The other pieces are \emph{outer}. As in \cite[Section~4]{wiegold}, call a segment of a $w^k$-circle \emph{consistent} if both of the arcs incident on the endpoints of the segment are oriented in the same direction. 

We will now double-count the total number of consistent segments. Firstly, as a consequence of Proposition \ref{proposition:nosourcesinkcycle}, every inner disc piece will have at least two consistent segments. So $2\chi(P)\leq sc(P)$ for every inner piece $P$, where $sc(P)$ denotes the number of consistent segments in the piece $P$. We hence have that the total number of consistent segments is at least $2 \sum_{\text{inner pieces}} \chi(P)$.

Secondly, note in a $w^k$-circle, any $|w|$ many arcs corresponding to consecutive junctures cannot be oriented in the same direction. Indeed, the orientations are given via the relative stacking $\rho$ and $x \cdot \rho(u)$ must take minimal and maximal values for exactly two of the prefixes. The orientations on the corresponding arcs will be in opposite directions (see ~\cite[Lemma $4.6$]{wiegold} for more details). Therefore, at least $2$ among any $|w|$ consecutive segments must be inconsistent, implying that we have at least $2\deg(S)$ inconsistent segments. Thus, we have at most $(|w|-2)\deg(S)$ many consistent segments.

Combining the two estimates for the number of consistent segments, we get $$\sum_{\text{inner pieces}} \chi(P) \leq \frac{1}{2} (|w|-2)\deg(S).$$

Finally, note that $S$ has exactly $\text{length}\of{u}$ many outer pieces and hence, we get the euler characteristic estimate, as follows.

\begin{align*}
&& -\chi(S) &=  \Fraction {1} / {2} (\text{length}\of{u} + \modulus{w} \deg(S)) - \sum_{\text{all pieces}} \chi(P) && \\
&&  &\geq \Fraction {1} / {2} (\text{length}\of{u}+ \modulus{w} \deg(S)) - \frac{1}{2} (|w|-2)\deg(S) - \sum_{\text{outer pieces}} \chi(P) && \\ && &\geq \Fraction {1} / {2} (\text{length}\of{u}+ \modulus{w} \deg(S)) - \frac{1}{2} (|w|-2)\deg(S) - \text{length}\of{u} && \\ && &= \deg(S) - \frac{1}{2} \text{length}\of{u}.
\end{align*}
\end{proof}

We briefly discuss different approaches. To prove asphericity, one can also use $w$-admisssible (planar) surfaces and the characteristic bound \cite[Theorem~D]{wiegold}. The difference between $w$-admissible surfaces and pictures is that for the former, one allows powers of $w$ and some $A$ and $B$ boundaries. Pictures can have any group element $g\in \normalClosure{w}$ as labelling the boundary. As Wilton recently mentioned to us, the power of stackings is not fully explored in the literature and it can be used to prove torsion-freeness of one-relator products of left-orderable groups. It would be interesting to further develop this machinery.


\section{Asphericity for One-Relator Products of Right-Orderable Groups}\label{section:asphericity}

In this section we construct an explicit classifying space for the one-relator product $A\ast B/\normalClosure{w}$, allowing us to prove a relative asphericity result (Theorem~\ref{lettertheorem:asphericity}) for the said one-relator product. 

Recall that $K$ was the cell complex obtained by gluing together a $K(A,1)$ and a $K(B,1)$ along the edge $e$. Fix a locally injective loop $l$ in the $1$-skeleton of $K$ representing the word $w\in \pi_1(X, \ast)$. Let $L$ be the complex obtained by attaching a $2$-cell $\sigma_w$ to $K$, whose attaching map is this loop $l$. We prove below that $L$ is aspherical which will be the key ingredient in our proof of Theorem \ref{corollary:torsionfree}.

\begin{theorem}
\label{theorem:asphericity}
Any continuous map from a 2-sphere $S$ to $L$, can be homotoped to another map whose image lies inside $K$, i.e. $\pi_2(L,K)=0$. Hence, $L$ is aspherical.
\end{theorem}

\begin{proof}
Let $\overline{f} \colon S \to L$ be a continuous map from a 2-sphere to $L$. First, via transversality, we homotope the map $\overline{f}$ so that the sphere $S$ inherits a natural cell structure encoding the map $\overline{f}$ (see \cite[Section 3]{FennRourke1996}, \cite[Section 3]{wiegold}, \cite{BuonchristianoRourkeSanderson_Transversality}). Make the map $\overline{f}$ transverse to a point $x$ in the interior of the open 2-cell $\sigma_w$ corresponding to $w$. The pre-image of a small neighborhood of $x$ will be a disjoin union of discs $D_1, D_2, D_3, ..., D_k$. After a homotopy, we may assume that the image $\overline{f}{D_i}$ is the entire 2-cell $\sigma_w$. We then get that $\overline{f}$ restricts to a map $f\colon S \setMinus \p{\Union_{i=1}^k D_i} \to K$. Without loss of generality, let $k$ be minimal over all maps in the homotopy class of $\overline{f}$. We will show $k=0$, proving that $\overline{f}$ can be homotoped to a map whose image lies in $K$.

Note that $(f, S)$ is a picture over the corresponding one-relator product. Suppose the picture is compressible, as witnessed by some pair of pants $P$. By compressing $P$, we obtain a picture $(f',S'=S)$ with $\overline{f} \simeq \overline{f'}$ (by Lemma \ref{homotopy}). Moreover, the preimage of a small neighborhood of $x$ is a disjoint union of $k-2$ discs contradicting the minimality of $k$. If the (oriented) picture $\Gamma := (f, S)$ is incompressible, it must contain a source/sink or elementary cycle by Proposition \ref{proposition:nosourcesinkcycle}. As $\Gamma$ is a picture on the sphere $S$, this contradicts Lemma \ref{elementarycycle}.

Therefore $\pi_2(L)=0$, and by Hurewicz and Whitehead, we must have that $L$ is aspherical. Indeed, let $p\colon \widetilde{L}\to L$ be the universal covering map, and let $\widetilde{K} := p^{-1}(K_A \sqcup K_B)$ be the pre-image under $p$ of the disjoint union of $K_A$ and $K_B$. Since the maps $K_A\to L$ and $K_B\to L$ are $\pi_1$-injective by \cite[~Theorem B]{wiegold}, $\widetilde{K}$ is a disjoint union of copies of the universal covers of $K_A$ and $K_B$. As $K_A$ and $K_B$ are classifying spaces, in particular they have contractible universal covers, and $\widetilde{L}\setminus \widetilde{K}$ consists only of $1$-cells and $2$-cells, the long exact sequence in homology for the pair $(\widetilde{L},\widetilde{K})$ gives that $H_i(\widetilde{L})=0$ for $i>2$. Finally, Hurewicz Theorem tells us $H_2(\widetilde{L})\cong \pi_2(L) = 0$, and by Whitehead we then have that $\widetilde{L}$ is contractible. Therefore, $L$ is aspherical.
\end{proof}

We now prove that for $w$ not a proper power, the one-relator product $A\ast B/\normalClosure{w}$ is torsion-free. For this, we require the following lemma from the theory of covering spaces (see~\cite{Hatcher2002, May1999}). We give its proof for completeness.
\newpage
\begin{lemma}\label{lemma:deck-transformation}
Let $X$ be a connected cell complex with universal covering map $p \colon \widetilde{X} \onto X$.
\begin{itemize}
\item Let $Y \subseteq X$ be a path-connected subspace, $\iota \colon Y \into X$ be the inclusion map, and $\iota_\ast \colon \pi_1\of{Y, y_0} \into \pi_1\of{X, x_0}$ be the induced map, where $y_0 \in Y$ and $x_0 := \iota\of{y_0}$. 
\item Let $Z_0$ be a path-connected component of $p^{-1}\of{Y}$ and $z_0 \in Z_0$ such that $p\of{z_0} = y_0$.
\end{itemize}
Consider the standard action of $G:= \pi_1(X, x_0)$ by deck transformations on $\widetilde{X}$. Let $g \in G$. Then, for any path-connected component $Z$ of $p^{-1}\of{Y}$, $g \cdot Z = Z$ if, and only if, $g \in h\iota_\ast\of{\pi_1\of{Y, y_0}}h^{-1}$  where $h$ is an element of $G$ satisfying $Z= h \cdot Z_0$.
\end{lemma}

\begin{proof}
Suppose $h \in G$ such that $Z = h \cdot Z_0$ and $g \in h \iota_\ast\of{\pi_1\of{Y, y_0}} h^{-1}$.  Then, $h^{-1} g h = [\alpha]$, where $\alpha$ is a loop based at $y_0 = x_0$ in $Y$. Lift $\alpha$ to a path $\beta$ in $Z_0$, such that $\beta\of{0} = z_0$.  The endpoint $\beta\of{1}$ of the path $\beta$ is $\p{h^{-1} g h} \cdot z_0$.  It follows that $g \cdot \p{h \cdot z_0} \in Z$, and therefore $Z \intersect g \cdot Z \neq \emptySet$. By connectedness, it follows $g \cdot Z = Z$.  

Conversely, suppose $g \cdot Z = Z$. There exists $h \in G$ such that $Z = h \cdot Z_0$. Then, $g \cdot \p{h \cdot z_0} \in Z =  h \cdot Z_0$. Let $\beta$ be the path from $z_0$ to $\p{h^{-1} g h} \cdot z_0$ in $Z_0$. Let $\alpha := p \circ \beta$, which yields a loop based at $y_0$ in $Y$. Then, $\alpha$ lifts to $\beta$, yielding $h^{-1} g h = [\alpha] \in  \iota^\ast\of{\pi_1\of{Y, y_0}}$. Thus, $g \in h \iota^\ast\of{\pi_1\of{Y, y_0}} h^{-1}$.
\end{proof}

\begin{corollary}
\label{corollary:torsionfree}
 Let $A, B$ be right-orderable groups and $w\in A\ast B$ be an element which is not conjugate into $A$ or $B$ and is not a proper power. Then, $A\ast B/\normalClosure{w}$ is torsion-free.
\end{corollary}

\begin{proof}
We follow a similar proof strategy to \cite[Theorem 2]{Forester_2005}. Let $\widetilde{L}$ be the universal cover of $L$ and let $L_A$ and $L_B$ be the pre-images of $K(A,1)$ and $K(B,1)$ with respect to the covering map. By Theorem~\ref{theorem:asphericity}, $\widetilde{L}$ is contractible. Again, note that any connected component of $L_A$ or $L_B$ (say $Z \subset L_A$) is contractible because the map $K(A,1)\to L$ is $\pi_1$-injective by \cite[~Theorem B]{wiegold}, so $Z$ is simply a copy of the universal cover of $K(A,1)$.

Contracting each of these disjoint sub-complexes to a point will give us a space $L'$ which is homotopy equivalent to $\widetilde{L}$, see ~\cite[Proposition $0.16$ and Exercise $0.27$]{Hatcher2002}. Thus, $L'$ is a contractible $2$-complex.

Suppose that $A*B/\normalClosure{w}$ has a non-trivial torsion element $g$. Now $g$ acts by deck transformations on $\widetilde{L}$ and $g\cdot L_A= L_A$ and $g\cdot L_B= L_B$. In particular, $g$ sends each connected component of $L_A$ or $L_B$ to another such component. Thus the action descends to an action on $L'$. The action on $L'$ is free away from the vertices obtained by collapsing the components of $L_A$ and $L_B$, and if $g\cdot v = v$ for such a vertex $v\in L'$, then we have $g\cdot Z = Z$ for some connected component of $L_A$ or $L_B$. By Lemma \ref{lemma:deck-transformation}, this means $g$ lies inside some conjugate of $A$ or $B$, which contradicts $A$, $B$ being torsion free. Applying this argument to the powers of $g$, we see that, the action of $\presentation{g}$ on $L'$ is free but $L'$ is a contractible $2$-complex, so this is impossible. Indeed, if the action of $F=\presentation{g}$ were free then $L^\prime/F$ would be a finite dimensional $K(F,1)$. However, finite cyclic groups have homology in infinitely many dimensions, so do not admit a finite dimensional classifying space.
\end{proof}

By using the $K\of{G, 1}$ space as constructed in Theorem~\ref{theorem:asphericity}, we are able to compute (co)-homology groups of $A\ast B/\normalClosure{w}$, and prove a relative version of the Lyndon Identity Theorem. The proofs for one-relator products of right-orderable groups are identical to their analogous results for one-relator products of locally indicable groups as proved by Howie in~\cite{Howie1984CohomologyOO}, and we therefore omit them.

\begin{corollary}
\label{corollary:cohomologycompute}
Let $A, B$ be right-orderable groups and $w\in A\ast B$ be an element which is not conjugate into $A$ or $B$ and is not a proper power. Let $G$ be the one-relator product $A\ast B/\normalClosure{w}$ and $\mathbf{1}$ denote the trivial subgroup of $G$. Then, for $k\geq3$,
\[H^{k}(G;\text{--}) = H^{k}(A;\text{--}) \oplus H^{k}(B;\text{--}), \]
 \[H_k(G;\text{--}) = H_k(A;\text{--}) \oplus H_k(B;\text{--}). \]

For $k = 1,2$, the (co)-homology groups of $G$ can be computed via the following exact sequences:
\[0\rightarrow H^{1}(G;\text{--}) \rightarrow H^{1}(A\ast B; \text{--})\rightarrow H^{0}(\mathbf{1}; \text{--})\rightarrow H^{2}(G; \text{--})\rightarrow H^{2}(A; \text{--})\oplus H^{2}(B; \text{--})\rightarrow 0,\]
\[0\rightarrow H_{2}(A; \text{--})\oplus H_{2}(B; \text{--}) \rightarrow H_{2}(G; \text{--})\rightarrow H_{0}(\mathbf{1}; \text{--})\rightarrow H_{1}(A\ast B; \text{--})\rightarrow H_{1}(G; \text{--})\rightarrow 0.\]
\end{corollary}

\begin{corollary}
\label{corollary:identitytheorem}
For right-orderable groups $A$ and $B$, let $N$ denote the normal closure of an alternating word $w\in A\ast B$ which is not a proper power. Let $G$ be the one-relator product $A\ast B/\normalClosure{w}$. Then, the abelianisation $N^{ab}$ of $N$ is isomorphic to $\integers G$ with generator $w\cdot [N,N]$.
\end{corollary}

\section{Classical Results for One-Relator Products of Right-Orderable Groups}\label{section:classical}

Let $A,B$ be right-orderable groups. We first prove that proper non-empty subwords of $w$ do not represent the identity element in $A*B/\normalClosure{w}$, generalising Weinbaum's subword theorem.

\begin{theorem}\label{theorem:subword}
Let $A, B$ be right-orderable groups and $w\in A\ast B$ be a word which is not conjugate into $A$ or $B$. Then, for any proper, non-trivial subword $u$ of $w$, we have $u \notin \normalClosure{w}$.
\end{theorem}

\begin{proof}

We give a proof by contradiction. To this end, suppose $u$ is a proper, non-trivial subword of $w$ such that $u \in \normalClosure{w}$, that is to say, $u = \prod_{i=1}^{N_u} g_iwg_i^{-1}$. Choose $u$ so that $N_u$ is minimal among all proper, non-empty subwords of $w$. Up to possibly taking a conjugate of $w$, one can assume without loss of generality that $u$ is a prefix of $w$. We notice that it is enough to prove the result when $w$ is not a proper power. Indeed, suppose $w=\hat{w}^k$ with $\hat{w}$ not a proper power. If $u$ is a power of $\hat{w}$, then by \cite[Corollary 6.1]{wiegold} $u\notin \normalClosure{\hat{w}^k}$. Otherwise, $u = \hat{w}^lu'$ with $u'$ being a proper, non-trivial subword of $\hat{w}$, so we can apply the result for $\hat{w}$ to get $u'\notin \hat{w}$, and hence $u\notin \normalClosure{\hat{w}^k}$.

Since $u\in \normalClosure{w}$, we can consider an oriented picture $\Gamma$, with the word $u$ on the boundary. Moreover, the representation $u = \prod_{i=1}^{N_u} g_iwg_i^{-1}$ gives a relative van Kampen diagram with $N_u$ many $2$-cells \cite[van Kampen's Lemma]{pictures}. The corresponding picture $\Gamma$, being the dual of this van Kampen diagram, will then have exactly $N_u$ many $w$-circles. In fact, we now note that one can also orient the arcs connected to the $\partial \Gamma$ using the relative stacking $\rho$. For any arc connecting a juncture corresponding to a prefix $u'$ of $u$ and a juncture corresponding to the prefix $u''$ of $w$, we orient it from $u'$ to $u''$ if $x \cdot \rho\of{u^\prime} > x \cdot \rho\of{u^{\prime\prime}}$. Notice that we cannot have $u'=u''$, so that we can always orient the arcs. Indeed, if $u'=u''$, then $v=u^{-1}w$ is a subword of $w$ such that $v\in \normalClosure{w}$ and $N_v<N_u$, contradicting the minimality of $N_u$. (See Figure \ref{fig:CompleteOrientation}).

\begin{figure}[ht]
    \centering
    \includegraphics[width=45mm]{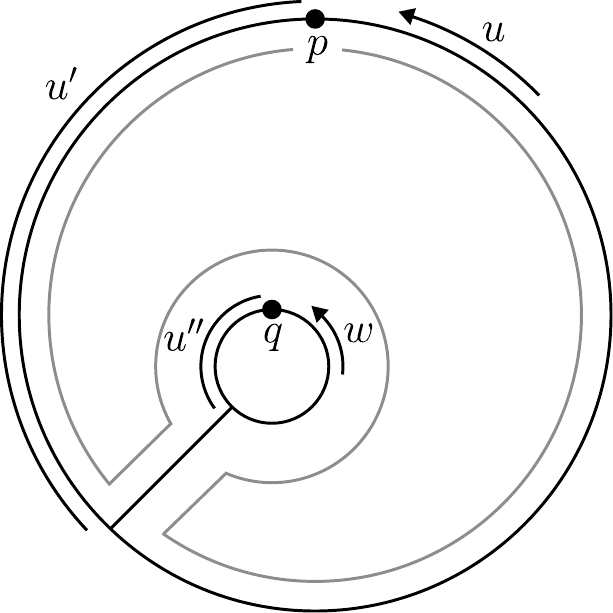}  
    \caption{\label{fig:CompleteOrientation}Following the grey path from the base point $p$, we read the word $u'(u'')^{-1}w^{-1}u''u'^{-1}u = v^{-1}$. Moreover, the path bounds a smaller number of $w$-boundary components.}
\end{figure}

Now that we have a picture with all arcs oriented, we will find a suitable oriented cycle yielding a contradiction. Notice that for every edge on the $u$-boundary connecting two junctures corresponding to prefixes $u'$ and $u''=u't$ we have $x \cdot \rho\of{u^{\prime\prime}}=x\cdot \rho\of{u^{\prime}t}$ except for the edge $e$, corresponding to the last letter of $u$.

\begin{figure}[ht]
    \centering
    \includegraphics[width=60mm]{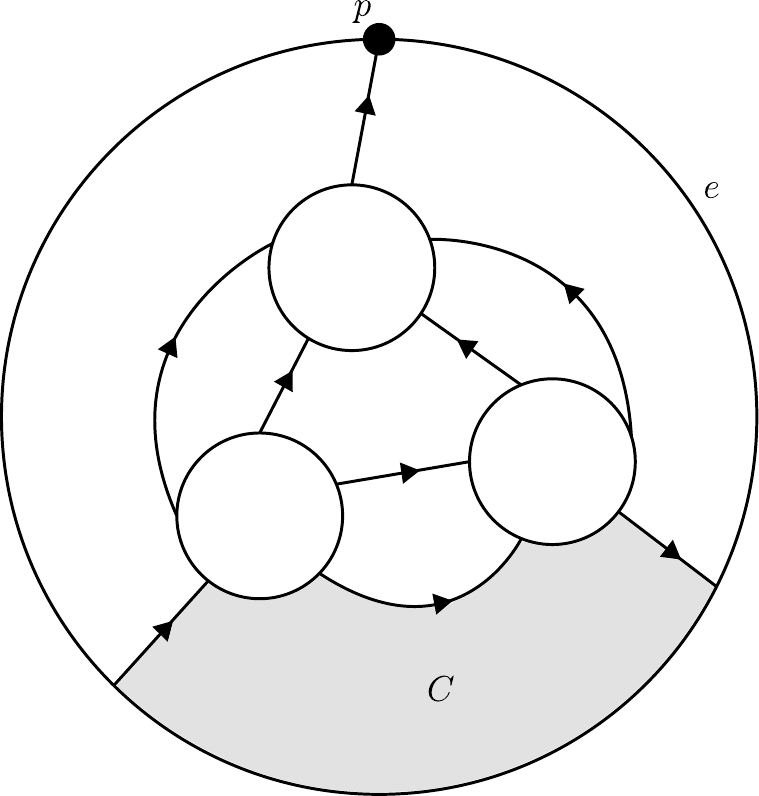}  
    \caption{\label{fig:SubwordThm}All arcs are oriented and we can find an elementary cycle $C$ which doesn't contain the edge $e$.}
\end{figure}

\begin{claim}
    There exists an elementary oriented cycle in $\Gamma$, which doesn't contain the edge $e$.
\end{claim}

\begin{proof}

Think of the picture $\Gamma$ as an oriented spherical graph with vertices being $\partial \Gamma$ and the $w$-circles, and edges being the arcs connecting them. We distinguish two cases.

\textbf{Case 1:} The vertex corresponding to $\partial \Gamma$ is a source or a sink.

Then, the corresponding graph will have one source/sink and hence we must have an oriented elementary cycle $C$ by Lemma \ref{elementarycycle}. If this cycle contains the edge $e$, then $\partial \Gamma$ must contain inwards and outgoing arcs contradicting that the corresponding vertex is a source/sink.

\textbf{Case 2:} The vertex corresponding to $\partial \Gamma$ is not a source or a sink.

Then $\Gamma$ tesselates the sphere and has no source or sink implying that we can find two elementary oriented cycles by Lemma \ref{elementarycycle}. One of these cycles must not contain the edge $e$, as there is only one elementary cycle containing $e$.
\end{proof}

Now that we have an oriented cycle $C$ each of whose arcs are oriented via the relative stacking $\rho$, we obtain our desired contradiction via Proposition \ref{proposition:nosourcesinkcycle}.
\end{proof}

We now prove a version of B. B. Newman's Spelling Theorem for one-relator products of right-orderable groups. 

\begin{theorem}
 \label{theorem:newmanspellingthm}
Let $A$, $B$ be right-orderable groups and consider the one-relator product $G=A*B/\normalClosure{w^k}$ where $w$ is not a proper power and not conjugate into $A$ or $B$. Let $u$ be a reduced alternating word in $A \ast B$ representing the identity element in $G$. Then, $u$ contains an alternating subword $v$ so that either $v$ or $v^{-1}$ is a subword of $w^k$ with 
\[
\emph{length}\of{v} > \Fraction k - 1 / {k} \emph{length}\of{w^k}.
\]
\end{theorem}

\begin{proof}
Since $u\in \normalClosure{w^k}$, we obtain a reduced orientated picture with some $w^k$-circles and such that $\partial \Gamma$ reads the alternating word $u$. We say that two arcs going from a fixed $w^k$-circle to $\partial \Gamma$ are \emph{consecutive} if they are consecutive both in $u$ and $w^k$ and they don't bound other $w^k$-circles. We will prove that one of the $w^k$-circles of our picture will have at least $(k - 1) \text{length}\of{w}+1$ many consecutive arcs connecting to $\partial \Gamma$. Then, taking $v$ to be the said part of the $u$-boundary proves the theorem.

Call a $w^k$-circle \emph{nice} if all the arcs connecting to $\partial \Gamma$ are consecutive (see Figure \ref{fig:NiceCircle}).

\begin{figure}[ht]
    \centering
    \includegraphics[width=100mm]{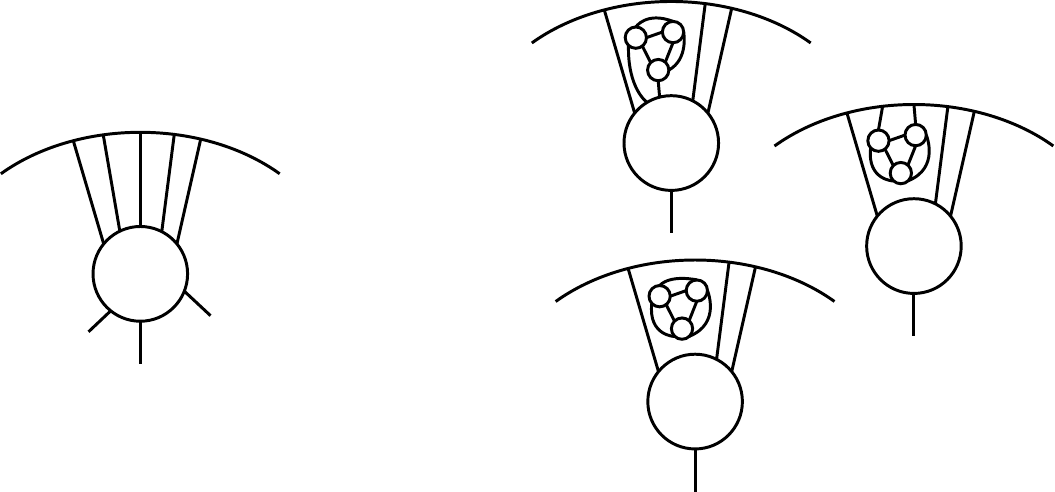}  
    \caption{\label{fig:NiceCircle}The left side of this figure shows a nice $w^k$-circle (the arc at the top is part of $\partial \Gamma$ and the three lines at the bottom connect to other $w^k$-circles not represented here). On the right of this figure, we can see three examples of $w^k$-circles which are not nice.}
\end{figure}

\begin{claim}
Let $l$ be a nice $w^k$-circle with fewer than $(k - 1) \text{length}\of{w}+1$ many consecutive arcs connecting to $\partial \Gamma$, then $l$ will admit both incoming and outgoing oriented arcs.
\end{claim}

\begin{proof}
Note that $l$ has at least $\text{length}\of{w}$ arcs connecting to other $w^k$-circles. Moreover, these arcs are consecutive in $l$. So it is possible to find a prefix of $w$ minimizing (respectively maximizing) the image of $x$ in a relative stacking for $w$. This prefix, by the definition of our orientation, will correspond to an incoming (respectively outgoing) arc.
\end{proof}

If all the $w^k$-circles are nice, then either we have found a $w^k$-circle with enough consecutive arcs or, by the claim above, each of the $w^k$-circles will admit both incoming and outgoing oriented arcs. For this second case, thinking of the picture $\Gamma$ as an oriented planar graph, by Lemma $\ref{elementarycycle}$, we must obtain an elementary cycle. This contradicts Proposition \ref{proposition:nosourcesinkcycle}.

Now, suppose the picture has some $w^k$-circle $l$ which is not nice so that it must admit two arcs to $\partial \Gamma$ which bound a region $R_l$ containing $n_l>0$ many $w^k$-circles (see Figure \ref{fig:Spelling}). Without loss of generality, assume that $n_l$ is minimal amongst all the $w^k$-circles that are not nice and all corresponding regions $R_l$ in the picture. If any of the $w^k$-circles inside $R_l$ is not nice, then the said boundary component will have fewer than $n_l$ $w^k$-circles inside a corresponding region contradicting the minimality of $n_l$. Therefore, the region $R_l$ must only contain nice $w^k$-circles.

\begin{figure}[ht]
    \centering
    \includegraphics[width=70mm]{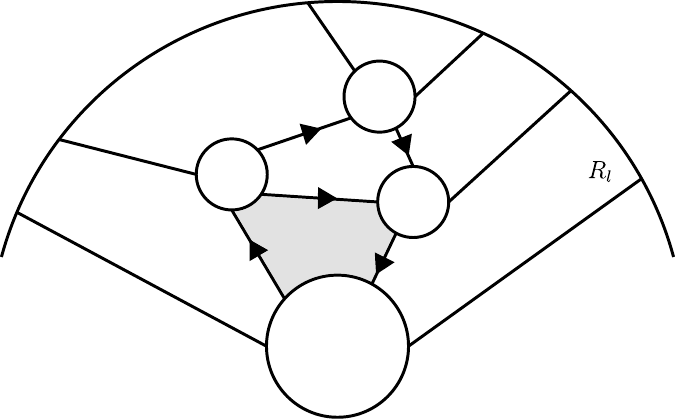}
    \caption{\label{fig:Spelling}All $w^k$-boundary components inside $R_l$ are nice, yielding an elementary oriented cycle.}
\end{figure}

Finally, think of the sub-picture corresponding to the region $R_l$ as an oriented spherical graph, with vertices being the $w^k$-circles and $\partial R_l$ and edges being the oriented arcs between them. By the claim above, every vertex corresponding to a $w^k$-circle must admit incoming and outgoing edges. Hence the corresponding graph has at most one source/sink, namely the vertex corresponding to $\partial R_l$. Again, by Lemma $\ref{elementarycycle}$, we must obtain an elementary cycle, contradicting Proposition \ref{proposition:nosourcesinkcycle}.
\end{proof}

\begin{remark}
A variation of the proof of Theorem \ref{theorem:newmanspellingthm} shows that you can get two disjoint alternating subwords $v_1$, $v_2$ of $u$ so that $v_i$ or $v_i^{-1}$ is a subword of $w^k$ and $\text{length}\of{v_i} > \Fraction k - 1 / {k} \text{length}\of{w^k}$.
\end{remark}

We recall that a group $G$ is \emph{hyperbolic relative} to a collection of subgroups $\set{H_1, \ldots, H_k}$ (in the sense of Osin) if $G$ is finitely presented with respect to $\set{H_1, \ldots, H_k}$ and the relative Dehn function of $G$ with respect to $\set{H_1, \ldots, H_k}$ is linear (see \cite[Definition $2.35$]{osin2004relatively}).

\begin{corollary}
\label{relativehyperbolicity}
 For right-orderable groups $A,B$, a one-relator product $A*B/\normalClosure{w^k}$ with torsion is hyperbolic relative to $\{A,B\}$.
\end{corollary}

\begin{proof}
 If the one-relator product has torsion, then $k\geq 2$, by Theorem \ref{corollary:torsionfree}. Theorem \ref{theorem:newmanspellingthm} implies that any word $u$ representing the identity element in $A*B/\normalClosure{w^k}$ must contain atleast half of $w^k$ or $w^{-k}$ as a subword. Therefore, $A*B/\normalClosure{w^k}$ has a linear relative Dehn function and is hence hyperbolic relative to $\set{A, B}$ (see for example \cite{osin2004relatively}).
\end{proof}

\begin{remark}
\label{Remark:relativehyperbolicity}
Corollary \ref{relativehyperbolicity} can also be proved without using Theorem \ref{theorem:newmanspellingthm} by adapting \cite[Theorem 4.1]{wiegold} to oriented pictures, as we shall now sketch.\\
Let $u$ be a word representing the identity in $A*B/\normalClosure{w}$. This means that $u\in \normalClosure{w}$, so there is a planar oriented picture $\Gamma$ with the word $u$ on the boundary, portraying $u$ as a product of conjugates of $w$. Assuming that this representation of $u$ as a product of conjugates is minimal, our picture must be incompressible. Note that we have $\chi\of{S} = 1-\text{area}\of{S}$ and $\deg\of{S} = k \cdot \text{area}\of{S}$. By Theorem \ref{theorem:characteristicestimateboundary}, we  obtain 
\[
-\chi\of{S} \geq \deg\of{S} - \Fraction 1 / 2 \text{length}\of{u}.
\]
Rewriting the above gives us 
\[
\text{area}\of{S} \leq \Fraction 1 / {2\p{k-1}} \text{length}\of{u} - \Fraction 1 / {k-1}. 
\]
Thus, $A \ast B/\normalClosure{w^k}$ satisfies a relative linear isoperimetric inequality and is hence hyperbolic relative to $\set{A, B}$ (see for example \cite{osin2004relatively}).
 \[\]
\end{remark}

\bibliographystyle{alpha}
\bibliography{References}

\end{document}